\documentclass[11pt]{elsarticle}
\usepackage[centertags]{amsmath}
\usepackage{amssymb}
\usepackage{amsthm}
\usepackage{amsfonts,amssymb}

\usepackage{mathrsfs}
\usepackage[nodots]{numcompress}

\usepackage{graphicx}

\theoremstyle{plain} 
\newtheorem{cor}{Corollary}

\newtheorem{pro}{Proposition}
\newtheorem{thm}{Theorem}
\newtheorem{conj}{Conjecture}

\theoremstyle{definition}

\theoremstyle{remark}

\DeclareMathOperator{\pd}{\partial}

\DeclareFontFamily{OT1}{pzc}{}
\DeclareFontShape{OT1}{pzc}{m}{it}{<-> s * [1.10] pzcmi7t}{}
\DeclareMathAlphabet{\mathpzc}{OT1}{pzc}{m}{it}

\DeclareMathSymbol{\R}{\mathalpha}{AMSb}{"52}
\DeclareMathSymbol{\C}{\mathalpha}{AMSb}{"43}
\newcommand{\mbb}[1]{\mathbb{#1}}

\newcommand{\N}{\mbb{N}}

\newcommand*{\Scale}[2][4]{\scalebox{#1}{$#2$}}%

\newcommand{\set}[1]{\left\{#1\right\}}
\newcommand{\comment}[1]{}

\newcommand{\bv}{\mathbf{v}}

\newcommand{\bw}{\mathbf{w}}

\newcommand{\beqn}{\begin{equation}}
\newcommand{\eeqn}{\end{equation}}
\newcommand{\balign}{\begin{align}}
\newcommand{\ealign}{\end{align}}
\newcommand{\bsube}{\begin{subequations}}
\newcommand{\esube}{\end{subequations}}

\begin{document}

\begin{frontmatter}



\title{Equivalence classes and Linearization of the Riccati and Abel chain}

\author{J.C. Ndogmo\corref{cor1}}
\ead{jean-claude.ndogmo@univen.ac.za}

\address{Department of Mathematics and Applied Mathematics\\
University of Venda\\
P/B X5050, Thohoyandou 0950, South Africa}

\author{Adrian M. Escobar-Ruiz}
\ead{escobarr@crm.umontreal.ca}

\address{
Centre de recherches math\'ematiques,
and D\'epartement de math\'ematiques \\
et de statistique, Universit\'e de Montreal,  C.P. 6128,
succ. Centre-ville,\\ Montr\'eal (QC) H3C 3J7, Canada}

\cortext[cor1]{Corresponding author}


\begin{abstract}
The problem of linearization  by point transformations is solved for equations in the generalized Riccati and Abel chain of order not exceeding the fourth. It is shown in particular that nonlinear third order and fourth order equations from the chain are not linearizable by any point transformations. The Lie pseudo-group of equivalence transformations for equations of arbitrary orders from the chain are then found, together with expressions for the transformed parameter functions. An important subgroup of the group of equivalence transformations found is considered and some associated equivalence classes are exhibited.
\end{abstract}


\end{frontmatter}

\section{Introduction}
\label{s:intro}
The sequences of nonlinear ordinary differential equations ({\sc ode}s)
\begin{subequations} \label{rabel}
\begin{align}
\left( \frac{d}{dx} + \lambda y \right)^k y &=0,\quad k \in \N := \set{1, 2, 3, \dots} \label{riccat1}      \\
\intertext{ and }
\left( \frac{d}{dx} + \lambda y^2 \right)^k y &=0,\quad k \in \N  \label{abel1}
\end{align}
\end{subequations}
where $\lambda$ is a given constant,  are commonly referred to as the Riccati and the Abel chains of equations, respectively. Indeed, the first order equations in \eqref{riccat1} and \eqref{abel1} are special cases of the standard Riccati and the standard Abel equations, respectively. These equations frequently occur in the theory of dynamical systems where they play an important role \cite{bruz1, bruz2, bruz3}.  One reason why the Riccati equation is also popular in physics and in many other mathematically based fields is that the similarity reduction of a second order linear {\sc ode} by scaling symmetries yields a first order Riccati equation. Higher order Riccati equations have also been shown to play a special role  in the construction of B\"{a}cklund transformations for integrable systems \cite{grund1}.\par

Over the recent years some Lie group related analysis of the Riccati and Abel chains have been performed in several papers. Thus in \cite{bruz-gdrias}, dynamical symmetries and related similarity reductions were obtained for equations \eqref{rabel} of low orders not exceeding the fourth. The third order Riccati equation in \eqref{riccat1} has also been  studied in \cite{carinena1} from a Lagrangian perspective, and their integrability has been investigated in \cite{carinena2} using the theory of Darboux and the extended Prelle-Singer methods.  Methods for identifying integrable cases of coupled Riccati and Abel chains \eqref{rabel} by means of nonlocal transformations have been proposed in \cite{lackshma1}.

Due to the simplicity of linear differential equations, and in particular their integrability and their many other familiar and nontrivial properties, a lot of interest has been devoted in the recent literature to the problem of linearization of differential equations, thus pursuing a work pioneered  by Lie \cite{lielin} and many others. The methods of linearization usually considered in such studies include the linearization by point or contact transformations, the generalized Sundman transformations as well as other methods (see \cite{pavel1Y18, soh&FM16, yulia2010iop, gungor17, ibra3, ibra4, DA.Lyakhov, muriel2010iop} and the references therein). However, a serious study of the linearization of the Riccati and Abel chain, especially by point transformations, does not seem to have been undertaken yet in the mathematical research literature.\par

Let $F=F(x,y)$ be a scalar function and consider the ordinary differential operator $\Omega=  \frac{d}{dx} + F(x,y)\equiv \Omega_F.$ This operator gives rise to iterative equations
\begin{equation} \label{gnlra}
 \Omega_F^n[y]\equiv \Omega^n[y] =0,
\end{equation}
obtained as iterations of the first order differential equation
\begin{subequations} \label{niter}
\begin{align}
\Omega[y] & \equiv  y_x+ F y =0,  \label{niter1}\\
\intertext{where $y = y (x),$}
\Omega^n [y] & = \Omega^{n-1}\left[\Omega [y] \right],\quad \text{ for $n\geq 1,$}  \label{niter2}
\end{align}
\end{subequations}
 and $\Omega^0 =I$ is the identity operator.  The function $F$ in \eqref{gnlra} is referred to as the parameter function of the differential operator $\Omega$. More recently some attempts have been made in \cite{chanpreserv} to find linearizing and order-preserving contact transformations as well as dynamical symmetries  for the generalized version of the Riccati and Abel chains of the form given by \eqref{gnlra}, which for simplicity  we shall still refer to  as the Riccati and Abel chain. Moreover, all these investigations were also applied to a coupled version of the Riccati and Abel chains, although the equations considered in the paper were of order not exceeding the third.\par

In this paper, we  solve the problem of linearization by point transformations of the Riccati and Abel chains \eqref{gnlra} for equations of order  $n$ such that $2 \leq  n \leq 4$.   These point transformations are more tractable than the nonlocal ones \cite{chanpreserv} as they preserve the symmetry algebra and immediately yield, in particular, the general solution of a given equation from that of its linearized counterpart.  For third order and fourth order equations from the general chain \eqref{gnlra}, we show in particular that they are not linearizable by point transformation unless they are already linear. We also state the latter result as a conjecture for all equations of order $n \geq 3,$ based on the fact that the iteration of a nonlinearizable equation is not likely to yield a linearizable one.\par

On the other hand we determine the Lie pseudo-group of equivalence transformations  for equations of all orders from the Riccati and Abel chain \eqref{gnlra} as well as the corresponding expression of the parameter function of the transformed equation under equivalence transformations. Canonical forms are also given for some of the most common types of parameter functions which are likely to be needed in Lie group classification problems. Also, a special subgroup of the Lie pseudo-group of equivalence transformations found is exhibited, and nontrivial equivalence classes of equations from the chain relative to this group are determined. Finally most of the results thus obtained are applied to special cases of the Riccati and Abel chains given by \eqref{riccat1} and \eqref{abel1} and new results about these simpler chains are derived.


\section{Linearization by point transformations}
\label{s:lin}

\subsection{Second order {\sc ode}s}
\comment{
Unlike nonlocal and other types of transformations, the linearization by point transformations of a differential equation has the interesting feature that it maps a linearizable equation to essentially the same equation in a different local coordinates system. As such, point transformations leave essentially invariant many important characteristics of the differential equations such as the symmetry algebra which is preserved up to isomorphism, and they yield a direct functional relation between the solution of the original equation and that of the transformed counterpart.\par
}
Algorithms for the linearization of {\sc ode}s of order not exceeding the fourth are well known \cite{lielin, ibra3, ibra4}.  Lie \cite{lielin} was indeed certainly the first to point out the fact that a scalar second order {\sc ode}
\begin{equation}\label{eq:nlodegn}
\Delta \equiv y_{xx}+ H(x, y, y_x)=0,
\end{equation}
linearizable by an invertible point transformation  must be at most cubic in the first derivative, that is, of the form
\begin{subequations} \label{lielincond}
\begin{align}
&\Delta \equiv y_{xx} + A y_{x}^{\,3}+ B y_{x}^{\,2}+ C y_x+
D =0 \label{lincond2},\\
\intertext{ for some functions $A, B, C$ and $D$ of $x$ and $y.$   When the necessary condition \eqref{lincond2} holds, the equation is linearizable if and only if
}
\Gamma_1 \left[ \Delta \right] & =0, \qquad \text{ and } \qquad  \Gamma_2\left[ \Delta \right]=0,  \label{lielincond1}\\
\intertext{where}
\begin{split}
\Gamma_1\left[ \Delta \right]&:=  3 A_{xx} -2 B_{xy} + C_{yy} - 3(C A)_x\\
    &\quad + 3(D A)_{y}+(B^2)_x + 3 A D_y - B C_y,
\end{split}\\
\begin{split}
    \Gamma_2\left[ \Delta \right]&:=  3 D_{yy}-2 C_{xy}+ B_{xx} - 3(D A)_{x} \\
    &\quad  +3(D B )_y- (C^2)_y - 3 DA_x+ C
    B_x.
\end{split}
\end{align}
\end{subequations}
Moreover, it also follows from Lie's linearization algorithm that a general point transformation
\begin{equation} \label{gen_ptrf}
z= \rho (x,y),\qquad w= \psi(x,y)
\end{equation}
maps the trivial equation $w''(z)=0$ to \eqref{lincond2} if and only if the functions $\rho$ and $\psi$ satisfy the determining equations which
are given for  $\rho_y= 0$  by

\begin{align} \label{deteqlin1}
\begin{split}
&\psi_{yy} = \psi_y B,\quad 2 \psi_{xy} = \rho_x^{-1} \psi_y \rho_{xx} + \psi_{y}C,\quad  \psi_{xx} = \rho_x^{-1} \psi_x  \rho_{xx}  + \psi_y D, \\
&2 \rho_x \psi_{xxx} -3 \rho_{xx}^2 - \rho_x^2 \left[ 4 (D_y + B D) - (2 C_x + C^2)\right]=0.\quad
\end{split}
\end{align}

On the other hand,  for $\rho_y \neq 0$ the determining equations are given by
\begin{align} \label{lielin2}
\begin{split}
 0  =\,&\rho_y \left(A \psi_x+\psi_{yy}   \right) -\psi_y \left(A \rho_x+\rho_{yy}\right)   \\
 \begin{split}
 0 =\,&\left(B \rho_y - A \rho_x\right) \left(\rho_y \psi_x-\rho_x \psi_y\right)-2 \rho_y \psi_y \rho_{xy}\\
 &+\left(\rho_x \psi_y -\rho_y \psi_x \right) \rho_{yy}+2 \rho_y^2 \psi_{xy}
 \end{split}\\
 0 =\,&\rho_y^2 \psi_{xx} -D \rho_y^2 \psi_y+\psi_x \left(C \rho_y^2-2 \rho_y \rho_{xy}+\rho_x \left(A \rho_x-B \rho_y+\rho_{yy}\right)\right) \\
 0 =\,&\rho_y^2 \rho_{xx} -D \rho_y^3+\rho_x \left(C \rho_y^2-2 \rho_y \rho_{xy}+\rho_x \left(A \rho_x-B \rho_y+\rho_{yy}\right)\right) \\
 \begin{split}
 0 =\,&2 \left(A B+A_y\right) \rho_x \rho_y+\left(B^2-4 A C+4 A_x-2 B_y\right) \rho_y^2\\
    &+6 A \rho_y \rho_{xy}-3 \left(A \rho_x+\rho_{yy}\right)^2+2 \rho_y \rho_{yyy}
 \end{split}\\
\begin{split}
 0 =\,& -3 A^2 \rho_x^3+3 \left(B^2-2 A C+2 A_x\right) \rho_x \rho_y^2-2 \left(3 A D-B_x+2 C_y\right) \rho_y^3\\
&+3 \rho_x \rho_{yy} \left(-2 B \rho_y+\rho_{yy}\right)+6 \rho_y \rho_{xy} \left(B \rho_y-2 \rho_{yy}\right) +6 \rho_y^2 \rho_{xyy}.
\end{split}
\end{split}
\end{align}
For second order equations from the chain \eqref{gnlra}, namely
\begin{equation} \label{2nd-o-ra}
0=\Omega^2[y]\equiv  y \left(F^2+F_x\right)+\left(2 F+y F_y\right) y_x+y^{(2)}
\end{equation}
where $F=F(x,y),$ Lie's conditions \eqref{lielincond1} reduce to
\begin{align*}
0& = 4 F_{yy} + yF_{yyy} \\
0&= 2(F- y F_y) F_{yy} + F_{xyy}.
\end{align*}
Solving this shows that $0 = \Omega^2[y]$ is linearizable if and only if
\begin{equation} \label{linzable2F}
F = \alpha_0(x)y + \beta(x)
\end{equation}
for some arbitrary functions $\alpha_0$ and $\beta$ of $x.$ \par

With this value of $F,$ and under the change of dependent variable $y = w e^{- \int_{z_0}^z \beta(t) dt},$ and by reusing the original notations $x$ and $y$ for variables in the transformed equation, the corresponding second order equation $\Omega_F^2 [y]=0$ from  \eqref{gnlra} is reduced to
\begin{equation}  \label{main2nd}
\alpha^2 y^3 + \alpha_x y^2 + 3 \alpha y y_x + y_{xx} =0,
\end{equation}
where $\alpha = \alpha_0 e^{- \int_{z_0}^z \beta(t) dt}.$ In fact \eqref{main2nd} is just $\Omega_F^2 [y] =0$ with $F = \alpha(x) y.$ This just says that second order equations in the chain \eqref{gnlra} with the function $F$ linear in $y$ and of the form \eqref{linzable2F} are equivalent under point transformations to an equation in the same chain with $F= \alpha(x) y.$ This type of equivalence between two given functions of the form $F=F(x,y),$ each labeling an equation from the chain \eqref{gnlra} can be readily established using the results on the equivalence transformations of \eqref{gnlra} obtained in Section \ref{s:eqvtrans}.\par

For the linearizable second order equation \eqref{main2nd}, the linearizing transformations of the form \eqref{gen_ptrf} are found after some calculations to be given as follows. For $\rho_y=0,$ the corresponding determining equations \eqref{deteqlin1} are incompatible and yield no solutions. For $\rho_y \neq0,$ solving the corresponding system \eqref{lielin2} yields a class of particular linearizing transformations given by

\begin{align} \label{linalpha}
\begin{split}
\rho (x,y) &=  \frac{x-k_5}{k_3\, y}- k_1 + T(x) \\
\begin{split}
\psi(x,y)  &= Q(x,y) \bigg[\big( (k_5-1) k_4 (k_5-x) + k_6 (x-1)\big) (x-k_5)  \\
&\quad - k_3 \left( (k_5-1) (k_4 k_1-k_2) (k_5-x) + k_6 k_1 (x-1)  \right) y \\
&\quad + k_3\, y \bigg[ \left( (k_5-1) k_4 (k_5-x) + k_6 (x-1)\right) T(x)\\
&\quad + (k_5-1) (k_5-x)\int \frac{k_6 \left(  k_1- T(x) \right)}{(k_5-x)^2} dx \bigg] \bigg],
\end{split}\\
\intertext{where}
T(s) &= \int \frac{(k_5-s) \alpha (s) }{ k_3   } \, ds \\
Q(x,y) &= \left[(k_5-1) k_3 (k_5-x) y \right]^{-1},
\end{split}
\end{align}

and the $k_j$ are arbitrary constants with $k_3 k_6 \neq 0$ and $k_5 \neq 1.$  We have thus established the following result.

\begin{thm} \label{t:lin2-ord}
Second order equations from the chain \eqref{gnlra}, namely equations of the form \eqref{2nd-o-ra}, where $F=F(x,y)$ is a given function, are linearizable by a point transformation if and only if $F= \alpha  y + \beta ,$ or equivalently $F= \alpha\, y,$ for some arbitrary functions $\alpha$ and $\beta$ of $x.$\par
Letting $F= \alpha\, y,$  a family of linearizing transformations for the corresponding  second order equations \eqref{main2nd} is given by \eqref{linalpha}.
\end{thm}

   Finding the symmetries of \eqref{main2nd} directly by the standard Lie approach through the solving of the determining equations can be a very intricate, if not impossible task. However, the linearizing  transformations \eqref{linalpha} or any particular case thereof can easily be used to obtain the full symmetry group of \eqref{main2nd} from those of its trivial counterpart. For finding these symmetries, we choose a particular case of \eqref{linalpha} which is as simple as possible by letting $k_3=k_6=1,$ and all the other constants $k_j$ equal to $0.$ Then considering the resulting linearizing transformations as a mere change of coordinates, and applying this to the symmetries of the trivial equation naturally yield the symmetries  of \eqref{main2nd} as a pull back of the local diffeomorphism defined by the change of coordinates. Two of the symmetries $\bv_j$ of \eqref{main2nd} obtained in this way, which we denote by $\bv_1$ and $\bv_2$ have  very simple expressions given by
\begin{align*}
\bv_1 & = x y \pd_x\, + \, y^2 (1 - x y \alpha) \pd_y  \\
\bv_2  &=   y(1-x) \pd_x\, + \, y^2 \left[ (x-1) y \alpha  -1 \right ] \pd_y.
\end{align*}
Apart from these two, all other symmetries of \eqref{main2nd} have lengthy expressions involving integrals of $\alpha(x).$ For instance, one of them which we denote by $\bv_3= \xi \pd_x + \phi \pd_y$ has expression
\begin{align*}
\xi &= \nu(x) \left[   (x-1) (x-y\, \mu(x)) + xy\, \nu(x)\right]\\
\begin{split}
\phi &=  - \frac{1}{x^2} \left[ (x-1)(x- y\, \mu(x) ) + x y\, \nu(x) \right] \\
      &\quad     \times \left[ x - y\, \mu(x) + x y\, (x y\, \alpha (x) -1) \nu(x)\right], \\
\end{split}
\intertext{where}
\mu(s) & = \int s \alpha(s) d s  \qquad \text{and}\qquad  \nu(s)  =  \int \frac{\mu(s)}{s^2} ds.\\
\end{align*}


\subsection{Higher order {\sc ode}s}
Based on similar ideas of Lie, Ibragimov, Meleshko, and collaborators recently extended Lie's algorithm \cite{lielin} to third order {\sc ode}s \cite{ibra3} and fourth order {\sc ode}s \cite{ibra4}. Although these algorithms often involve lengthy calculations, their implementation is straightforward.
It follows from such algorithms \cite{ibra3} that the linearization conditions for the third order equation $\Omega^3[y] =0$ reduce to
\begin{align*}
0&=4 F_y + y F_{yy}\\
0&=3F_y + y F_{yy}\\
0&= y F_y^2 -2 y^2 F_y F_{yy} + 3 F (F_y + y F_{yy}) + 3 (F_{xy} + y F_{xyy}).
\end{align*}
This shows that $\Omega^3[y]=0$ is linearizable if and only if  $F=F(x),$ in which case the equation itself is already linear.  Indeed,  for $F=F(x),$  $\Omega_F$ is a linear operator, and thus all equations in the chain $\Omega_F^n[y]=0$ are linear.\par

Similarly it follows from \cite{ibra4} that the linearization conditions for the fourth order equation $\Omega^4[y]=0$ include
\begin{align*}
0&=5 F_y + y F_{yy},\\
0&= F_{yy}.
\end{align*}
This shows that $F= F(x)$ must hold. Substitution of the latter equality in the other remaining determining equations for linearization shows that $F= F(x)$ is also a sufficient condition of linearization on $F(x,y).$ Since for $ F=F(x)$ the equations  $\Omega^n[y] = 0$ are linear for all $n \geq 1,$ we have just found that as for the third order equation, for given $F=F(x,y)$ the fourth order equation $\Omega_F^4[y]=0$ is linearizable only when it is already a linear one.\par

The above discussion thus shows that $\Omega_F^n[y]=0$ is linearizable for $n=3,4$ if and only if $F_y = 0.$ One question we wish to answer is whether this statement also holds for all $n >4.$ This boils down to finding out if for some given function $F= F(x,y)$ and for some $n\geq 4$ the equation $\Omega_F^n[y]=0$ may not be linearizable while $\Omega_F^{n+1}[y]=0$ is. However, from the structure of an iterative equation, it is very unlikely that the iteration of a non-linearizable equation will yield a linearizable one. This leads us to suggest the following result.
\begin{conj} \label{conj1}
For $n \geq 3,$ the equation
\[
\Omega^n[y] = 0
\]
is linearizable if and only if  $F_y =0,$ in which case it is already a linear equation.
\end{conj}
It is worth recalling that based on existing linearization algorithms, we have proved this conjecture for the orders $n=3 \text{ or } 4.$ \par

 It should be noted that when $F_y=0,$ that is, for linear equations from the chain, the change of variable $y = w e^{- \int_{x_0}^x F(t) dt}$ transforms  $\Omega_F^n [y]=0$ into the trivial equation $w^{(n)}=0,$ for all $n \geq 1.$ In fact for all arbitrary functions $F=F(x,y)$ the nonlocal change of variable
\begin{equation}\label{nloc-chge}
y = w e^{- \frac{1}{2}\int_{x_0}^x F(t, y(t)) dt}
\end{equation}
also maps any equation $\Omega^n [y]=0$  in the chain  to the trivial equation $w^{(n)}=0.$ The drawback with the nonlocal transformation \eqref{nloc-chge} is however that it seems a bit too strong, in the sense that it doesn't preserve important features of the original equation such as its symmetry algebra, and doesn't also provide a  practical relationship between the solution of a given equation and that of its transformed version.

\section{Equivalence transformations}
\label{s:eqvtrans}

In this section we are interested in finding the most general invertible point transformations of the form
\begin{equation} \label{gen_ptrfv2}
x= \rho (z , w),\qquad y= \psi(z,w)
\end{equation}
where $\rho$ and $\psi$ are some  functions satisfying $\rho_z \psi_w - \rho_w \psi_z \neq 0$ and   the form of the equation $\Omega_F^n [y]=0$ is preserved for a given value of $n,$ and up to the arbitrary function $F.$ In other words, the functions $\rho$ and $\psi$ are to be found in such a way that the set of all equations $\Omega_F^n [y]=0,$  where $F$ takes on arbitrary values, is invariant under \eqref{gen_ptrfv2}. Equivalently, each equation $\Omega_F^n [y]=0$  for a given  parameter function $F$ should be transformed under \eqref{gen_ptrfv2} into an equation $\Omega_Q^n [w]=0$ of another chain with new parameter function $Q$ whose expression is also to be found.\par

\subsection{Determination of equivalence transformations}
The procedure we adopt here for finding equivalence transformations is a direct one. We subject the equation $\Omega_F^n [y]=0$ to a point transformation \eqref{gen_ptrfv2} and request that the resulting equation be invariant in the above described sense. The details of this procedure for the third order equation $\Omega_F^3 [y]=0$ are as follows. First we note that   this third order equation  has explicit expression
\begin{equation}
\begin{split}
\label{gracO3}
&0= y^{(3)} +y_x \left(3 \left(F^2+F_x\right)+y \left(3 F F_y+2 F_{xy}\right)\right) \\
&\quad\; + y \left(F^3+3 F F_x+F_{xx}\right) +y_x^2 \left(3 F_y+y F_{yy}\right)+\left(3 F+y F_y\right) y_{xx}.
\end{split}
\end{equation}
When \eqref{gracO3} is subjected to \eqref{gen_ptrfv2}, the new term $-3 \rho_w (w_{zz})^2$ appears in the transformed equation. The vanishing of this term shows that one must have $\rho= \rho(z).$  With this new expression for $\rho,$ another new term that appears in the transformed version of the equation is $3 \psi_{ww} w_z w_{zz} / \psi_w.$ The vanishing of this term also shows that $\psi (z,w) = R(z) w + S(z),$ for some arbitrary functions $R \neq 0$ and $S$ of $z.$ With these values of $\rho$ and $\psi,$ the coefficient of $w_{zz}$ in the transformed equation has expression
\[  3 F \rho_z+(S+R\, w) F_{S+R\, w} \rho_z+\frac{3 R_z}{R}-\frac{3 \rho_{zz}}{\rho_z}.
\]
Comparing this expression with the coefficient of $y_{xx}$ in the original equation \eqref{gracO3} shows that the parameter of the transformed equation must be of the form
\[ Q=  F( \rho, S+R\, w ) \rho_z+ \frac{R_z}{R}-\frac{\rho_{zz}}{\rho_z}.  \]

Comparing again the coefficients of $w_{zz}$ in the iterative equation $\Omega_{Q}[w]=0$ and in the transformed version of \eqref{gracO3} shows that the function $S= S(z)$ must identically vanish. Then the difference between the resulting equations, that is, the iterative equation $\Omega_{Q}^3[w]=0$ and the transformed version of  $\Omega_F^3 [y]=0$ has expression
\[
\Scale[1.0]{
\begin{split}
&\frac{w_z \left(3 \rho_{zz}^2-2 \rho_z \rho_{zzz}\right)}{\rho_z^2} -\frac{w \left[6 R \rho_{zz}^3-3 \rho_z \rho_{zz} \left(R_z \rho_{zz}+2 R \rho_{zzz}\right)\right]}{R \rho_z^3} \\
&+ \frac{F R
\rho_z^2 \left(-3 \rho_{zz}^2+2 \rho_z \rho_{zzz}\right)+\rho_z^2 \left(2 R_z \rho_{zzz}+R \rho_{zzzz}\right)}{R \rho_z^3}.
\end{split}
}
\]

For this expression to vanish identically, it is necessary and sufficient that the coefficient of $w_z$ in the expression identically vanishes, and this is achieved for  $\rho$ given by
\begin{equation}\label{rho}
\rho (z) = \frac{k_2 z + k_1}{k_4 z + k_3}, \qquad {k_4 =0 \text{ or } 1,}
\end{equation}
where the $k_j$ are arbitrary constants for $j=1,2,3.$  With these expressions for  $\rho$ and $\psi= R(z) w,$ one sees that \eqref{gen_ptrfv2} transforms \eqref{gracO3} into another iterative equation $\Omega_{Q}^3[w]=0$ with parameter function
\begin{equation}\label{newF}
Q (z,w) = F(\rho, R\, w ) \rho_z+ \frac{R_z}{R}- \frac{\rho_{zz}}{\rho_z},
\end{equation}
where $\rho$ is given by \eqref{rho} and as before $R$ is an arbitrary function of $z.$ The determination of the equivalence transformations and the transformed parameter for the orders $n=2$ and $n>3$ proceeds exactly along the same lines of reasoning as for the order $n=3.$\par

 For the order $n=1,$ it turns out that every transformation of the general form \eqref{gen_ptrfv2} maps the first order iterative equation $\Omega_F [y]=0$ to $\Omega_{Q_1}[w]=0,$ with
\begin{equation}
\label{Transfparn=1} Q_1= \frac{F(\rho, \psi) \psi \rho_z + \psi_z}{w\left[ F(\rho, \psi) \psi \rho_w + \psi_w \right]}.
\end{equation}
We have thus obtained the following result.
\begin{thm} \label{th:eqv-trfo}\mbox{$\quad$}
\begin{enumerate}

\itemsep=1mm
\item[(a)] For $n =1,$ the equivalence transformations of \eqref{gnlra} are given by any transformation of the form \eqref{gen_ptrfv2} and the parameter function  $Q_1$ of the transformed equation is given by \eqref{Transfparn=1}.

\item[(b)] For $n\geq 2,$ the equivalence transformations of \eqref{gnlra} are given by \eqref{gen_ptrfv2}, with
\begin{subequations} \label{gen-eqv}
\begin{align}
\rho &= \frac{k_2 z + k_1}{k_4 z + k_3} \qquad (k_4 =0 \text{ or } 1),\qquad \psi =  R\, w \label{gen-eqv1}\\
\intertext{ and the parameter function $Q_n (z,w)$ of the transformed equation has expression}
Q_n (z,w)&=  F(\rho, R\, w ) \rho_z+ \frac{R_z}{R}- \left( \frac{n-1}{2} \right)\frac{\rho_{zz}}{\rho_z}.  \label{gen-eqv2}
\end{align}
\end{subequations}
In \eqref{gen-eqv1}, $k_1, k_2,$ and $k_3$ are arbitrary constants and $R= R(z)$ is an arbitrary function, while the expression of $\rho$ in \eqref{gen-eqv2} is given by \eqref{gen-eqv1}.
\end{enumerate}
\end{thm}
\subsection{Canonical form of some parameter functions}\label{s:canonic}

Given that each equation in a Riccati and Abel chain can be uniquely identified with its parameter function, the equivalence transformations of Riccati and Abel chains naturally induce and equivalence relation on the set of all functions $F=F(x,y)$ considered as parameter functions of equations the form \eqref{gnlra}. Thanks to Theorem \ref{th:eqv-trfo}, it is now possible to reduce the parameter function $F$ of a given equation from the Riccati and Abel chain \eqref{gnlra} to a much simpler parameter function corresponding to an equivalent equation. Such a simplified parameter function may be referred to as the canonical form of all parameter functions in an equivalence class of {\sc ode}s belonging to Riccati and Abel chains. The knowledge of these canonical forms is usually crucial in the Lie group classification of equations, and this is the primary motivation for which we are interested here in them. We shall therefore restrict our attention to equations of order $n \geq 2,$ as the symmetry group of first order {\sc ode}s is generally of no interest.
\begin{cor} \label{cor:1}
Assume $n \geq 2$ in \eqref{gnlra}. Let $F_0(x,y)$ be a fixed function, and let $\varepsilon \neq 0, a, b, c,$ and $d$ be arbitrary constants.
\begin{enumerate}
\itemsep=1.5mm
\item[(a)]  The function $F(x,y) = \varepsilon F_0 (a x + b, c y + d)$ is equivalent to the function $Q_n(z,w) = F_0 (b + \frac{a z}{\varepsilon}, d+ w)$ for $c \neq 0,$ and to $Q_n(z,w)=0$ for $c=0.$
\item[(b)]  The function $F(x,y) = F_0 (x,y) + A(x) y + B(x),$ where $A(x)$ and $B(x)$ are given functions  is equivalent to the function
    \[ Q_n(z,w) =  F_0 (z, \varepsilon \exp {\textstyle (\int_{z_0}^z - B(t)dt)}   w) +   \varepsilon \exp {\textstyle (\int_{z_0}^z - B(t)dt)} A(z) w.\]
\item[(c)]  The function $F(x,y) = F_0 (x,A(x) y),$ where $A(x)$ is a given nonzero function  is equivalent to the function  $Q_n(z,w)=  F_0 (z,w) - A_z(z)/A(z).$
\end{enumerate}

\end{cor}
\begin{proof}
With the notations of Theorem \ref{th:eqv-trfo},   the functions $Q_n(z,w)$ are obtained as the transformed parameter function of $\Omega_F^n [y]=0$ under \eqref{gen-eqv1} and given by \eqref{gen-eqv2}.
For Part (a),  the precise expression for $Q_n (z,w)$ is obtained for $c\neq0$ by letting  $R(z) = \frac{1}{c} (\varepsilon \rho_z)^{\frac{n-1}{2}}$ and $\rho(z) = z/ \varepsilon.$  For $c=0,$ we may assume that $F=F(x).$ In this case letting $\rho (z)=z$ and $R(z) = \exp {\textstyle (\int_{z_0}^z - F(t)dt)}$ yields $Q_n(z,w)=0.$ Similarly, for part (b), $Q_n (z,w)$ is obtained from \eqref{gen-eqv2} with $\rho(z) = z$ and $R(z) = \exp (  \int_{z_0}^z  - B(t) dt ).$ In Part (c)  the function $Q_n(z,w)$ is obtained from \eqref{gen-eqv2} by letting $\rho(z)=z$ and $R(z) = 1/ A(z).$ This completes the proof of the corollary.
\end{proof}

\subsection{Equivalence subgroup preserving functions $F= \alpha(x) y$}  \label{ss:eqvsub}
It has been established in Section \ref{s:lin} that a second order equation \eqref{gnlra} labeled by a parameter function $F=F(x,y),$ that is, an equation of the explicit form \eqref{2nd-o-ra} is linearizable by point transformations if and only if $F = \alpha y + \beta$ for some arbitrary  functions $\alpha$ and $\beta$ of $x.$ Moreover, any such equation \eqref{2nd-o-ra} is equivalent to one in which $F= \alpha y.$ It becomes therefore natural to examine the equivalences classes of higher order equations \eqref{gnlra} labeled by parameter functions of the form $F = \alpha y.$\par
Let us denote by $G_1$ the Lie pseudo-group of  point transformations of the most general form \eqref{gen_ptrfv2}  acting on the space of independent and dependent variables. Let $G_0$ denote the subgroup of $G_1$ consisting of point transformations of the form \eqref{gen-eqv1}. Hence by construction $G_0$ is the largest group of point transformations leaving invariant  \eqref{gnlra}. Finally, let $G$ be the subgroup of $G_0$ leaving invariant equations \eqref{gnlra} with parameter functions $F = \alpha y.$ In other words, $G$ transforms any equation $\Omega_F^n[y]=0$ with $F = \alpha y$ to another equation $\Omega_{Q_n}^n[y]=0$ with $Q_n= \delta y,$ where as usual $\alpha$ and $\delta$ are functions of the independent variable $x.$\par
It follows from Conjecture \ref{conj1} that contrary to the case of second order equations, for $n>2$ none of the equations $\Omega_F^n[y]=0$ with $F = \alpha y,$ and at least for $n=3$ or $4$ is equivalent under $G_1$ to the corresponding equation in the chain with $F=0,$ that is, to the trivial equation. Moreover, it also turns out that under $G_1$ equations in the chain \eqref{gnlra} with $F = \alpha y$ do not form a single equivalence class for $n \geq 3,$ contrary to what has been established in Theorem \ref{t:lin2-ord} for second order equations. Indeed, let $\theta, u,v, r,$ and $s$ be constants such that $ \theta \neq0,\;u \neq 0,\;\text{ and } r \neq s.$ If we denote by $L_{n,F}$ the symmetry algebra of $\Omega_F^n[y]=0,$ then it is immediately found for $n=3$ or $4$ that one has for instance

\begin{align} \label{dimLnF}
\begin{split}
\dim L_{n, F} &= 3\; \text{ if } F= \theta y, \quad  \dim L_{n, F} = 1\; \text{ if } F= (u x + v) y,\\
\dim L_{n, F} &= 0\; \text{ if } F= \frac{x+r}{x+ s}y.
 \end{split}
\end{align}
The distinct values of $\dim L_{n, F}$ thus obtained show that a classification of functions $F= \alpha y$ under $G$ should yield meaningful inequivalent classes.

\begin{pro}\label{p:eqv-grp-G}
The subgroup $G$ of $G_1$ made up of equivalence transformations of \eqref{gnlra} where the parameter functions are of the form  $F = \alpha(x) y$ consists of transformations of the form \eqref{gen_ptrfv2}, in which
\begin{subequations}\label{subgroupG}
\begin{align}
\rho &= \frac{k_2 z + k_1}{k_4 z + k_3} \quad (k_4 =0 \text{ or } 1),\quad \text{and}\quad\psi =  r_0 \rho_z^{(n-1)/2} \, w, \label{subgroupG1}\\
\intertext{ where $r_0$ is an arbitrary nonzero constant. The corresponding expression of the transformed parameter function is given by}
Q_n (z,w) &= r_0\,  \alpha (\rho)\,  \rho_z^{(n+1)/2}w. \label{subgroupG2}
\end{align}
\end{subequations}

   \end{pro}

\begin{proof}
 It follows from \eqref{gen-eqv2} that under \eqref{gen-eqv1}, a function of the form $F= \alpha y$ is transformed into
   \begin{equation} \label{trfo-alphaf}
   Q_n (z,w)=  R\, w \alpha (\rho) \rho_z   + \frac{R_z}{R}- \frac{(n-1)\rho_{zz}}{2\, \rho_z}.
   \end{equation}
Therefore  $F= \alpha y$ is transformed into a function of the same form $Q_n (z,w)= \delta (z)\, w(z)$ if and only if
\begin{align*}
0 &= \frac{R_z}{R}- \frac{(n-1)\rho_{zz}}{2\, \rho_z},
\end{align*}
that is if and only if $R= r_0 \rho_z^{(n-1)/2}$  for some constant of integration $r_0.$ Substituting the latter expression for $R$ into \eqref{gen-eqv1} and \eqref{trfo-alphaf} shows that the subgroup $G$ of $G_0$  as well as the corresponding expression of the transformed  parameter function $Q_n(z,w)$ are indeed given by \eqref{subgroupG}.
\end{proof}

Let $\gamma \neq 0$ and $k_3, k_4$  be constants such that $(k_3, k_4)\neq (0,0).$ Consider the following sets of functions of the form $F= \alpha y$ in which $n\geq 2$ is as usual the order of the corresponding equation labeled by $F$ in \eqref{gnlra}.

\begin{subequations}
\begin{align}
E_1  &= \set{\frac{\gamma}{(k_4 x+ k_3)^{\,n+1}}\,y\; \colon \quad ( \gamma \neq 0 ) } \\
E_2   &=\set{   \frac{\gamma (p x + q)}{(k_4 x+ k_3)^{\,n+2}}\,y\; \colon \quad p k_3- k_4 q \neq 0,\; p \neq 0     } \\
E_3   &=\set{  \frac{\gamma}{(p x + q)(k_4 x+ k_3)^{\,n}}\,y\; \colon \quad p k_3- k_4 q \neq 0,\; p \neq 0      }  \\
E_4   &=\set{  \frac{\gamma (x + q)}{(x + s)(k_4 x+ k_3)^{\,n+1}}\,y\; \colon \quad  q\neq s,\; k_4 s - k_3\neq 0  }.
\end{align}
\end{subequations}
\begin{pro} \label{p:eqvclass}
Let $\theta, u,v, r,$ and $s$ be constants such that
\[
\theta \neq0,\quad u \neq 0,\quad \text{ and } r \neq s.
\]

Under the action of the group $G$ of equivalence transformations defined by \eqref{subgroupG1}, the following are,  in each case,  equivalence classes of the indicated functions.
\begin{enumerate}
\item[(a)] $E_1,$ for functions $F= \alpha_1  y,$ with $\alpha_1(x)= \theta.$

\item[(b)] $E_2,$ for functions $F= \alpha_2 y,$ with $\alpha_2(x)= (u x+ v).$

\item[(c)] $E_3,$ for functions $F= \alpha_3 y,$ with $\alpha_3(x)=   1/(u x+ v).$
\item[(d)] $E_4,$ for functions $F= \alpha_4 y,$ with $\alpha_4(x)=    ( x+ r)/ (x+ s).$
\end{enumerate}
Moreover the equivalence classes $E_j$ are pairwise nonequivalent under the action of $G.$
\end{pro}
\begin{proof}
It follows from Proposition \ref{p:eqv-grp-G} that each group element $g$ of $G$ can be identified with a pair $(r_0, \rho)$ given by \eqref{subgroupG1} and where it is of course also assumed that $\Delta \equiv k_2 k_3 - k_1 k_4\neq 0.$ It also follows from \eqref{subgroupG2} that the action of $g$ on a function $F= \alpha(x)y$ can be reduced to an action of $g$ on the corresponding function $\alpha = \alpha (x)$  and this may thus be denoted by
\begin{subequations} \label{actionG}
\begin{align}
g\cdot \alpha &= r_0\,  \alpha (\rho)\,  \rho_x^{(n+1)/2} =  \alpha (\rho)\, \frac{\omega}{(k_4 x+ k_3)^{n+1}}, \label{actionG1} \\
\intertext{using the fact that}
\rho_z^{(n+1)/2} & = \Delta^{(n+1)/2} /(k_4 x+ k_3)^{n+1},    \label{actionG2}
\end{align}
\end{subequations}
and thus $\omega = r_0 \Delta^{(n+1)/2}$ is an arbitrary constant.   It is therefore straightforward that $g\cdot \alpha_j$ is of the form prescribed by $E_j$ for each $j=1,\dots, 4.$ Conversely, it is also true that for each parameter function of the form $F =\beta_j y \in E_j,$ there is some $g\in G$ such that $g\cdot \alpha_j = \beta_j.$ We verify this for $j=3,$ noting that the other cases can be verified in a similar way. One has
\begin{align*}
g\cdot \alpha_3 (x) &= (r_0, \rho) \cdot \frac{1}{(u x+ v)} =   \alpha_3 (\rho) \frac{\omega}{(k_4 x+ k_3)^{n+1}} \\
  &= \frac{1}{(k_2 u + k_4 v) x + (k_1 u + k_3 v) } \frac{\omega}{(k_4 x+ k_3)^{n}}.
\end{align*}
Thus given $\alpha_3$ and $\beta_3= \gamma\, (p x + q)^{-1}(k_4 x+ k_3)^{\,-n},$ and the arbitrariness of $\omega,$ in order to find $g$ satisfying $g \cdot \alpha_3 = \beta_3,$ it suffices to find $\rho$ of the given form \eqref{subgroupG1}  such that
\[
k_2 u + k_4 v = p,\qquad \text{ and }\quad   k_1 u + k_3 v = q.
\]
But since $u\neq 0,$ such a linear fractional transformation is easily found up to two free parameters $k_1$ and $k_3.$ This shows that $E_3$ is exactly the equivalence class of $F= \alpha_3 y.$ \par
 On the other hand, it clearly follows from  the forms of the functions $F = \beta_j y$ in each set $E_j$ as well as the specified restrictions on the constants these  functions involve, that the sets $E_j$ are pairwise nonequivalent, and this completes the proof of the proposition.
\end{proof}
Although the equivalence transformations given by the group $G_0$ and even more by its subgroup $G$ are quite weak, especially due to their limited effect on the dependent variable $y,$ it clearly follows from \eqref{dimLnF} that the inequivalent classes $E_1, E_2,$ and $E_4$ remain inequivalent even under the action of the most general point transformations \eqref{gen_ptrfv2}. In other words, for functions $F= \alpha y,$ in the absence of the Lie group classification of equations \eqref{gnlra} or its  group of equivalence transformations for arbitrary orders $n,$ Proposition \ref{p:eqvclass} shows  that the subgroup $G$ can be used to identify important nontrivial equivalence classes for this particular family of equations.\par

\section{Application to ordinary Riccati and Abel chains}
\label{s:applicat}

  The special cases of the generalized Riccati and Abel chain \eqref{gnlra} given by \eqref{riccat1} and \eqref{abel1} are much better known and have been studied as already mentioned in several papers (see for instance  \cite{grund1, bruz-gdrias, carinena1, carinena2, chanpreserv} and the references therein.) To distinguish these special cases from the Riccati and Abel chain \eqref{gnlra}, we shall refer to \eqref{riccat1} and  \eqref{abel1} as the ordinary Riccati and the   ordinary Abel chains, respectively. In this section we briefly review these chains  in the light of some of the results obtained thus far in this paper. \par

  When the constant parameter $\lambda$ in \eqref{riccat1} and \eqref{abel1} equals zero, the Riccati and Abel chains are clearly identical and reduce to sequences of trivial equations. We therefore assume henceforth that $\lambda \neq 0.$ It then follows from Theorem \ref{th:eqv-trfo} that the change of  variable $x= z / \lambda $ reduces each of  \eqref{riccat1} and \eqref{abel1} to the same chain in which $\lambda=1.$ In order words, up to a point transformation one may always assume that [?? the] parameter function $F(x,y)$ in \eqref{gnlra} satisfies $F(x,y)= y$  in \eqref{riccat1} and $F(x,y)=y^2$ in \eqref{abel1}. Moreover, each equation of order $n$ in the Riccati chain \eqref{riccat1} is equivalent under point transformation to an equation in another Riccati chain with parameter function $Q_n(z,w)$ of the general form
  \begin{align*}
  Q_n(z,w) &= \rho' R w + \frac{R'}{R} - \frac{n-1}{2} \frac{\rho''}{\rho'}.  \\
  \intertext{ For the Abel chain \eqref{abel1}, a similar fact holds and the corresponding transformed parameter function takes the form}
  Q_n(z,w) &= \rho' R^2 w^2  + \frac{R'}{R} - \frac{n-1}{2} \frac{\rho''}{\rho'},
  \end{align*}
where the functions  $\rho= \rho(z)$ and $R= R(z)$ are those specified in Theorem \ref{th:eqv-trfo}. If however, we are only interested in equivalence transformations of the form \eqref{gen-eqv} preserving functions $F= \alpha y,$ for arbitrary functions $\alpha= \alpha(x),$ then it follows from Proposition \ref{p:eqvclass} that the Riccati chain \eqref{riccat1} can be viewed as the canonical form of chains corresponding to the class $E_1$ described by Proposition \ref{p:eqvclass}.\par

   On the other hand, it follows from Theorem \ref{t:lin2-ord} and Conjecture \ref{conj1} (proven for third and fourth order equations) that for orders $n \geq 2,$ none of the equations in the Abel chain \eqref{abel1} is linearizable by point transformation, while for the Riccati chain \eqref{riccat1}, only the second order equation is linearizable. Assuming as usual that $\lambda=1$ in \eqref{riccat1}, it follows from \eqref{linalpha} that the corresponding linearizing transformation in this case is given by
   \begin{equation*}
   z= \frac{x}{y} - \frac{x^2}{2}, \qquad w = \frac{x-1}{y} - \frac{x(x-2)}{2}.
   \end{equation*}
  These results show in particular that,  contrary to what has often been stated in the recent literature, at least for the orders $n$ such that $2 \leq n \leq 4,$ none of the equations from the chains \eqref{riccat1} and \eqref{abel1} is of maximal symmetry, except for the second order equation in \eqref{riccat1}.  In fact by direct calculation up to the order ten, the general symmetry vectors $\bv_n$ and $\bw_n$ for $n$th order equations from \eqref{riccat1} and \eqref{abel1} respectively, are given  by
\begin{alignat*}{2}
\bv_n &= \left[ k_1 + x(k_2+ k_3 x)\right]\pd_x + \left[ n k_3 - y (k_2 + 2 k_3 x)\right]\pd_y,& \qquad &(n \geq 3)  \\
\bw_n &=  \left( k_1 + k_2 x\right)\pd_x - (k_2\, y /2 ) \pd_y,& \qquad &(n \geq 2)
\end{alignat*}
where the $k_j$ are arbitrary constants. In particular,  $\bw_n$ does not depend on $n$ and the low dimensions of the symmetry algebras generated by $\bv_n$ and $\bw_n$ confirm the fact that the corresponding equations cannot be linearized by point transformations.\par

Nevertheless, it is easy to see, as is well known, that each $n$th order equation from the Riccati chain can be obtained as the Lie reduction under the scaling symmetry of the trivial $(n+1)$th order equation $w^{(n+1)}=0.$ Such a reduction is given explicit by the sequence of transformations $w= e^u,$ and then $y= \dot{u}\equiv \pd u / \pd x,$ and shows that the ordinary Riccati chain \eqref{riccat1} consists of exactly integrable equations.


\section{Concluding Remarks}
\label{s:conclusion}

One fact to be outlined at this point concerns the linearization problem considered in this paper.  This problem applies to equations from the generalized Riccati and Abel chain \eqref{gnlra} of order $n\geq 2.$ The problem is completely solved in this paper for equations of order $n$ with $2 \leq n \leq 4,$ and an answer is provided in Conjecture \ref{conj1} for equations of order $n \geq 5.$ No attempt has however been made to investigate the linearization of first order equations from the chain \eqref{gnlra} because as is well known \cite{DA.Lyakhov}, not only all first order {\sc ode}s which can be solved for the derivative are linearizable by a point transformation, but also finding such a transformation is generally as hard as solving the equation itself.\par

On the other hand, the solution to the linearization problem for second  order equations from the chain has pointed out a new problem that is worthwhile considering. Namely, the determination of the group of equivalence transformations for equations from the  chain corresponding to the more specific case where $F= \alpha(x) y$ in \eqref{gnlra}. It has indeed been  shown in this paper that second order equations from the chain are linearizable by a point transformation if and only if the equation is equivalent to one for which the function $F$ has the latter expression, and a transformation mapping the resulting linearizable equation to its trivial counterpart was obtained in \eqref{linalpha}. It thus seems interesting to see what the equivalence group for this specific value of $F$ looks like for equations of order three and four, and even for second order ones. The problem of finding the said group of equivalence transformations is just a generalization of that of finding the group $G$ found in Section \ref{ss:eqvsub}, not as a subgroup of the group $G_0$ anymore,  but  of the full Lie pseudo-group $G_1$ given by the most general point transformation \eqref{gen_ptrfv2}.  It should also be noted that $n$th order equations from the chain \eqref{gnlra} with $F= \alpha(x) y$ are in fact equivalent under \eqref{gen-eqv1} to equations with $F$ of the more general form  given by the function $Q_n(z,w)$ in \eqref{trfo-alphaf}.

\label{s:declare}

\section*{Acknowledgement}
\label{s:acknowledge}
JCN acknowledges financial support from the NRF Incentive Funding for
Rated Researchers [Grant Number 97822], and from the University of Venda
[Grant Number I538].

\bibliographystyle{model1-num-names}

\end{document}